\documentclass[12pt]{article}
\usepackage{a4}
\usepackage{epsfig}
\usepackage{amsfonts}
\usepackage{amsthm}
\newtheorem{theorem}{Theorem}
\newtheorem{lemma}[theorem]{Lemma}
\newtheorem{problem}{Problem}
\begin{document}
\title{Coloring plane graphs with independent crossings}
\author{Daniel Kr\'al'\thanks{Institute for Theoretical Computer Science (ITI), Faculty of Mathematics and Physics, Charles University, Malostransk\'e n\'am\v est\'\i{} 25, 118 00 Prague 1, Czech Republic. E-mail: {\tt kral@kam.mff.cuni.cz}. Institute for Theoretical computer science is supported as project 1M0545 by Czech Ministry of Education.}\and
        Ladislav Stacho\thanks{Department of Mathematics, Simon Fraser University, 8888 University Dr, Burnaby, BC, V5A 1S6, Canada. E-mail: {\tt lstacho@sfu.ca}. This research was supported by NSERC grant 611368.}}
\date{}
\maketitle
\begin{abstract}
We show that every plane graph with maximum face size four whose
all faces of size four are vertex-disjoint is cyclically $5$-colorable.
This answers a question of Albertson whether graphs drawn in the plane
with all crossings independent are $5$-colorable.
\end{abstract}

\section{Introduction}

Coloring of graphs embedded in surfaces, in the plane in particular, attracts a lot of attention
of researchers in graph theory. The famous Four Color Theorem~\cite{bib-appel76+,bib-robertson97+}
asserts that every graph that can be drawn in the plane with no crossings is $4$-colorable.
It is natural to ask what number of colors is needed to color graphs that can be embedded
in the plane with a restricted number of crossings. If every edge is crossed by at most
one edge (such graphs are called $1$-embeddable and we restrict our attention solely to such graphs throughout this paper), Ringel~\cite{bib-ringel65}
conjectured that six colors suffice. This conjecture was answered in affirmative
by Borodin~\cite{bib-borodin84,bib-borodin95}.

Albertson~\cite{bib-albertson08} considered graphs with even more restricted structure of crossings. 
Two distinct crossings are independent if the end-vertices of the crossed pair of edges are mutually
different. In particular, if all crossings are independent, then each edge is crossed by at most one other
edge. Albertson showed that every graph drawn in the plane with at most 3 crossings is $5$-colorable (note
that the complete graph of order five can be drawn in the plane with a single crossing) and conjectured~\cite{bib-albertson08,bib-albertson} that
every graph that can be drawn in the plane with all its crossings independent is $5$-colorable.
In this paper, we prove his conjecture.

The coloring problem that we study is closely related to the notion of cyclic coloring. A coloring
of vertices of an embedded graph is {\em cyclic} if any two vertices incident with the same face
receive distinct colors. Let us show how the original problem can be expressed using this notion.
Let $G$ be a plane graph with all its crossings independent. We can assume (by adding edges
if necessary) that all faces of $G$ that do not contain a crossing have size three and
those that contain a crossing have size four. Remove now all edges that are crossed
by another edge. Clearly, a cyclic coloring of the obtained graph $G'$ is a proper coloring
of the original graph $G$ and vice versa. The assumption that all crossings of $G$ are independent
translates to the fact that all faces of $G'$ with size four are vertex disjoint. Hence, our main
result can be stated as follows:

\begin{theorem}
\label{thm-main}
Let $G$ be a plane graph with faces of size three and four only.
If all the faces of size four are vertex-disjoint, then $G$ is
cyclically $5$-colorable.
\end{theorem}

Before we proceed with proving Theorem~\ref{thm-main}, let us survey known results on cyclic colorings of
plane graphs. Since the maximum face size is a lower bound on the number of colors needed in a cyclic
coloring, it is natural to study the number of colors needed to cyclically color a plane graph
as a function of its maximum face size $\Delta^*$. If $\Delta^*=3$, then the graph is a triangulation and
the optimal number of colors is four by the Four Color Theorem. If $\Delta^*=4$, then the optimal
number of colors six by results of Borodin~\cite{bib-borodin84,bib-borodin95}; the optimality is witnessed
by the prism over $K_3$. For larger values of $\Delta^*$, the Cyclic Coloring Conjecture of Ore and
Plummer~\cite{bib-ore69+} asserts that the optimal number of colors is equal
to $\left\lfloor 3\Delta^*/2\right\rfloor$ (the optimality is witnessed by a drawing of $K_4$
with subdivided edges). After a series of papers~\cite{bib-borodin92,bib-borodin99+} on this problem,
the best general bound of $\left\lceil 5\Delta^*/3\right\rceil$
has been obtained by Sanders and Zhao~\cite{bib-sanders01+}.
Amini, Esperet and van den Heuvel~\cite{bib-amini08+} cleverly used a result
by Havet, van den Heuvel, McDiarmid and Reed~\cite{bib-havet07+,bib-havet08+}
on coloring squares of planar graphs and showed that
the Cyclic Coloring Conjecture is asymptotically true in the following sense: for every $\varepsilon>0$,
there exists $\Delta_\varepsilon$ such that every plane graph of maximum face size
$\Delta^*\ge\Delta_\varepsilon$ admits a cyclic coloring with at most
$\left(\frac{3}{2}+\varepsilon\right)\Delta^*$ colors.

There are two other conjectures related to the Cyclic Coloring Conjecture of Ore and Plummer.
A conjecture of Plummer and Toft~\cite{bib-plummer87+} asserts that every $3$-connected
plane graph is cyclically $(\Delta^*+2)$-colorable. This conjecture is known to be true
for $\Delta^*\in\{3,4\}$ and $\Delta^*\ge 18$, see~\cite{bib-enomoto01+,bib-hornak99+,bib-hornak00+,bib-hornak+}.
The restriction of the problems to plane graphs with a bounded maximum face size is removed
in the Facial Coloring Conjecture~\cite{bib-kral05+} that asserts that vertices of every plane
graph can be colored with at most $3\ell+1$ colors in such a way that every two vertices
joined by a facial walk of length at most $\ell$ receive distinct colors. This conjecture
would imply the Cyclic Coloring Conjecture for odd values of $\Delta^*$. Partial results
towards proving this conjecture can be found in~\cite{bib-havet+,bib-havet++,bib-kral05+,bib-kral07+}.

\section{Preliminaries}
\label{sect-notation}

The proof of Theorem~\ref{thm-main} is divided into several steps. We first identify
configurations that cannot appear in a counterexample with the smallest number
of vertices. Later, using a discharging argument, we show that a plane graph avoiding
all these configurations cannot exist. In particular,
vertices and faces of a counterexample are assigned
charge whose total sum is negative and which is redistributed preserving its total sum.
Lemmas~\ref{lm-f55}--\ref{lm-v6} claim that the final amount of charge of every vertex and
every face is non-negative which excludes the existence of a counterexample and yields a proof
of Theorem~\ref{thm-main}.

We now introduce notation used throughout the paper. Let us start
with some general notation. A vertex of degree $d$ is referred to as a {\em $d$-vertex} and
a face of size $d$ as a {\em $d$-face}. A {\em cyclic neighbor} of a vertex 
$v$ is a vertex lying on the same face as $v$ and the {\em cyclic degree} of $v$
is the number of its cyclic neighbors.

Our goal is to prove
Theorem~\ref{thm-main}. We assume that the statement of the theorem is false and consider
a counterexample with the smallest number of vertices; such a counterexample
is referred to as {\em minimal}, i.e., a minimal counterexample $G$ is a plane
graph with faces of size three and four such that all $4$-faces of $G$
are vertex-disjoint, $G$ has no cyclic $5$-coloring and any graph $G'$
satisfying assumptions of Theorem~\ref{thm-main} with a smaller number
of vertices than $G$ has a cyclic $5$-coloring.

A vertex $v$ of a minimal counterexample $G$ is {\em pentagonal} if
the degree of $v$ is five, $v$ is incident with no $4$-face and every
neighbor of $v$ is incident with a $4$-face. A $4$-face incident
with a neighbor of a pentagonal vertex $v$ is said to be {\em close} to $v$
if it contains an edge between two consecutive neighbors of $v$;
a $4$-face incident with a neighbor of a pentagonal vertex that is not close
is {\em distant}. If $f$ is close/distant to a vertex $v$, then we also
say that $v$ is close/distant to $f$.
A pentagonal vertex is {\em solitary} if no $4$-face is close to it. 

\begin{figure}
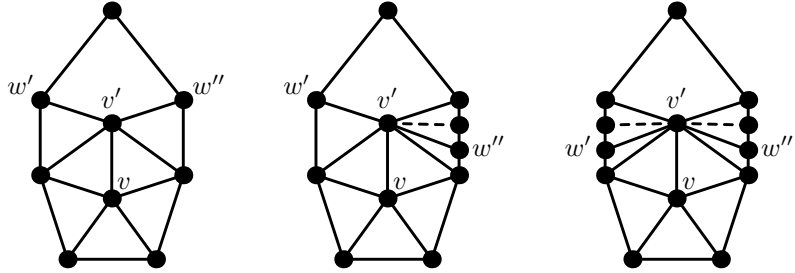

\begin{center}
\epsfbox{plan-ind.21}
\hskip 6mm
\epsfbox{plan-ind.22}
\hskip 6mm
\epsfbox{plan-ind.23}
\end{center}
\caption{Examples of a pentagonal vertex $v$ adjacent to a vertex $v'$ of
         degree five, a one-sided vertex and a double-sided vertex (in this order).}
\label{fig-sided}
\end{figure}

Let $v$ be a pentagonal vertex and $v'$ a neighbor of it. Let $w'$ and $w''$
be the common neighbors of $v'$ and another neighbor of $v$ (see Figure~\ref{fig-sided}).
If the $4$-face incident with $v'$ contains both $w'$ and $w''$, then
the degree of $v'$ is five. If the $4$-face contains one of the vertices
$w'$ and $w''$, then $v'$ is said to be {\em one-sided}, and if the $4$-face
incident with $v'$ contains neither $w'$ and $w''$, then $v'$ is {\em double-sided}.
Observe that if a pentagonal vertex is adjacent to a vertex of degree five,
it must also be adjacent to a double-sided vertex (otherwise, some of
the $4$-faces incident with its neighbors would not be vertex-disjoint).

\section{Reducible configurations}
\label{sect-reduce}

In this section, we show that a minimal counterexample cannot contain
certain substructures which we refer to as {\em configurations}.
Let us start with the following simple observation.

\begin{lemma}
\label{lm-sep}
A minimal counterexample $G$ does not contain a separating cycle of length two or
three.
\end{lemma}

\begin{proof}
Assume that $G$ contains a separating cycle $C$ of length two or three.
Let $G'$ and $G''$ be the subgraphs lying in the interior and the exterior of
the cycle $C$ (including the cycle $C$ itself). If $C$ is of length two, remove one
of the two parallel edges bounding $C$ from $G'$ and $G''$. By the minimality
of $G$, both $G'$ and $G''$ have a cyclic $5$-coloring. The colorings
of $G'$ and $G''$ readily combine to a cyclic $5$-coloring of $G$.
\end{proof}

\noindent We will use Lemma~\ref{lm-sep}
as follows: if we identify some vertices of a minimal counterexample,
Lemma~\ref{lm-sep} guarantees that the resulting graph is loopless
as long as every pair of the identified vertices have a common neighbor.
Indeed, if a loop appeared, the two identified vertices with their common
neighbor would form a separating cycle of length three.

We next show that the minimum degree of a minimal counterexample
is at least five.

\begin{lemma}
\label{lm-mindeg}
A minimal counterexample $G$ does not contain a vertex $v$ of degree four or less.
\end{lemma}

\begin{figure}
\begin{center}
\epsfbox{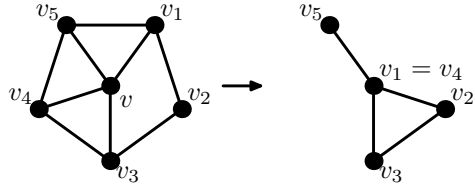}
\end{center}
\caption{A vertex of degree four with five cyclic neighbors and its reduction.}
\label{fig-mindeg}
\end{figure}

\begin{proof}
If the cyclic degree of $v$ is less than five, let $G'$ be the graph
obtained by removing $v$ from $G$ and triangulating the new face.
By the minimality of $G$, $G'$ has a cyclic $5$-coloring.
Since the cyclic degree of $v$ is less than five, this coloring
can be extended to a cyclic $5$-coloring of the original graph $G$.
Hence, we can assume that the cyclic degree of $v$ is five. In particular,
the degree of $v$ is four and $v$ is contained
in a $4$-face (see Figure~\ref{fig-mindeg}).

Let $v_1,\ldots,v_5$ be the neighbors of $v$. By symmetry we can assume that
the $4$-face incident with $v$ is $vv_1v_2v_3$. Let $G'$ be the graph
obtained from $G$ by removing the vertex $v$ and identifying the vertices
$v_1$ and $v_4$ to a new vertex $w$, see Figure~\ref{fig-mindeg}.
Note that the vertex $w$ is contained in at most one $4$-face
since the $4$-face incident with $v_1$ becomes a $3$-face in $G'$.
Since the vertices $v_1$ and $v_4$ have a common neighbor,
the graph $G'$ is loopless by Lemma~\ref{lm-sep}.

By the minimality of $G$, $G'$ has a cyclic $5$-coloring.
Since two of the neighbors of $v$ (the vertices $v_1$ and $v_4$)
are assigned the same color and the cyclic degree of $v$ is five,
the coloring can be extended to a cyclic $5$-coloring of $G$.
\end{proof}

Our next step is to show that all vertices of degree five that appear in a minimal counterexample 
must be pentagonal or incident with a $4$-face.

\begin{lemma}
\label{lm-pentagonal}
Every vertex $v$ of degree five in a minimal counterexample $G$
is either pentagonal or incident with a $4$-face.
\end{lemma}

\begin{figure}
\begin{center}
\epsfbox{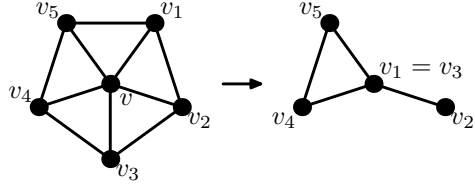}
\end{center}
\caption{A non-pentagonal vertex of degree five incident with no $4$-face and its reduction.}
\label{fig-pentagonal}
\end{figure}

\begin{proof}
We proceed as in the proof of Lemma~\ref{lm-mindeg}. Consider a $5$-vertex $v$
incident with $3$-faces only such that one of its neighbors is not incident
with a $4$-face. Let $v_1,\ldots,v_5$ be the neighbors of $v$ and $v_1$
a neighbor not incident with a $4$-face. Remove $v$ and identify
vertices $v_1$ and $v_3$ (see Figure~\ref{fig-pentagonal}).
Since the vertex $v_1$ is not incident
with a $4$-face in $G$, the new vertex is contained in at most one $4$-face.
By the minimality of $G$, the new graph can be cyclically $5$-colored and
this coloring readily yields a coloring of $G$.
\end{proof}

In the next lemma, we show that no $4$-face of a minimal counterexample contains two adjacent
vertices of degree five.

\begin{lemma}
\label{lm-square}
A minimal counterexample $G$ does not contain a $4$-face with two adjacent
vertices of degree five.
\end{lemma}

\begin{figure}
\begin{center}
\epsfbox{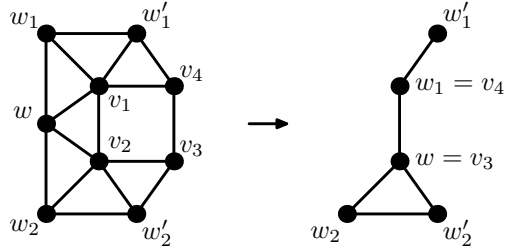}
\end{center}
\caption{A reduction of a $4$-face with two adjacent $5$-vertices.}
\label{fig-square}
\end{figure}

\begin{proof}
Assume that $G$ contains a $4$-face $v_1v_2v_3v_4$ such that the degrees of 
$v_1$ and $v_2$ are five. Let $w$ be the common neighbor of $v_1$ and $v_2$,
$w_1$ and $w'_1$ the other neighbors of $v_1$ (named in such a way that $w'_1$
is a neighbor of $v_4$) and $w_2$ and $w'_2$ the other neighbors of $v_2$.
See Figure~\ref{fig-square}.

Let $G'$ be the graph obtained by removing the vertices $v_1$ and $v_2$ and
identifying the vertices $w$ and $v_3$ and the vertices $w_1$ and $v_4$.
Clearly, the graph $G'$ is loopless (as the graph $G$ has no separating
$3$-cycles by Lemma~\ref{lm-sep}) and all its $4$-faces are vertex-disjoint.

By the minimality of $G$, $G'$ has a cyclic $5$-coloring. Assign the vertices
of $G$ the colors of their counterparts in $G'$. Next, color the vertex $v_2$:
observe that two of its $6$ cyclic neighbors have the same color and one
is uncolored. Hence, $v_2$ can be colored. Since the vertex $v_1$ has
$6$ cyclic neighbors and two pairs of its cyclic neighbors have the same
color, the coloring can also be extended to $v_1$.
\end{proof}

In the next two lemmas, we show that a $4$-face of a minimal counterexample
cannot contain a vertex of degree at most six adjacent to a close pentagonal vertex.

\begin{lemma}
\label{lm-pentagon5}
A minimal counterexample $G$ does not contain a vertex of degree five contained
in a $4$-face that is adjacent to a close pentagonal vertex.
\end{lemma}

\begin{figure}
\begin{center}
\epsfbox{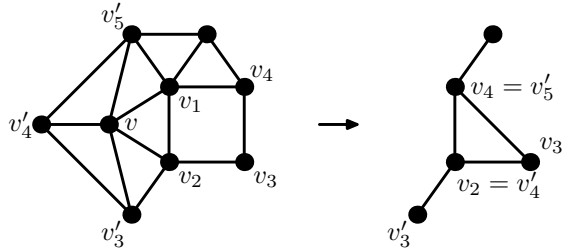}
\end{center}
\caption{A reduction of a $4$-face with a vertex of degree five adjacent to a close pentagonal vertex.}
\label{fig-pentagon5}
\end{figure}

\begin{proof}
Assume that $G$ contains a $4$-face $v_1v_2v_3v_4$ such that $v_1$ has degree
five and is adjacent to a close pentagonal vertex $v$. Let $v_1,v_2,v'_3,v'_4,v'_5$
be the neighbors of $v$ (see Figure~\ref{fig-pentagon5}).
Let $G'$ be the graph obtained by removing the vertices $v$ and $v_1$
and identifying the vertices $v_2$ and $v'_4$ and the vertices $v_4$ and $v'_5$.
Since every pair of identified vertices has a common neighbor, $G'$
is loopless by Lemma~\ref{lm-sep}. The $4$-faces of $G'$ are also
vertex-disjoint.

By the minimality of $G$, the graph $G'$ has a cyclic $5$-coloring. Assign
the vertices of $G$ the colors of their counterparts in $G'$. We next color
the vertex $v_1$ with an available color (the cyclic degree of $v_1$ is six,
it has a pair of neighbors colored with the same color and an uncolored neighbor)
and then the vertex $v$ (its cyclic degree is five and it has a pair of neighbors 
colored with the same color). The existence of this coloring contradicts that $G$ is a counterexample.
\end{proof}

\begin{lemma}
\label{lm-pentagon6}
A minimal counterexample does not contain a vertex of degree six contained
in a $4$-face that is adjacent to a close pentagonal vertex.
\end{lemma}

\begin{figure}
\begin{center}
\epsfbox{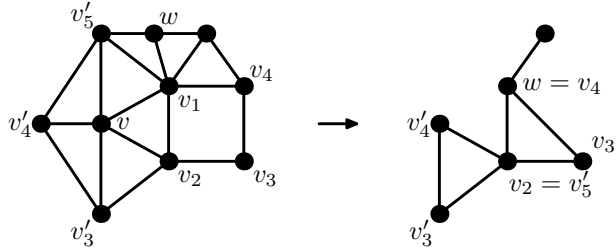}
\end{center}
\caption{A reduction of a $4$-face with a vertex of degree six adjacent to a close pentagonal vertex.}
\label{fig-pentagon6}
\end{figure}

\begin{proof}
Assume that $G$ contains a $4$-face $v_1v_2v_3v_4$ such that $v_1$ has degree
six and is adjacent to a close pentagonal vertex $v$. Let $v_1,v_2,v'_3,v'_4,v'_5$
be the neighbors of $v$ and $w$ the common neighbor of $v_1$ and $v'_5$ (since all $4$-faces are vertex disjoint, both faces containing the edge $v_1v'_5$ have size three and the vertex $w$ must exist).
Also see Figure~\ref{fig-pentagon6}.
Let $G'$ be the graph obtained from $G$ by removing the vertices $v$ and $v_1$
and identifying the vertices $v_2$ and $v'_5$ and the vertices $v_4$ and $w$.
Since every pair of identified vertices has a common neighbor, $G'$
is loopless by Lemma~\ref{lm-sep}. The $4$-faces of $G'$ are also
vertex-disjoint.

By the minimality of $G$, the graph $G'$ has a cyclic $5$-coloring. Assign
the vertices of $G$ the colors of their counterparts in $G'$. We next color
the vertex $v_1$ with an available color (the cyclic degree of $v_1$ is seven,
it has two pairs of neighbors colored with the same color and an uncolored neighbor)
and then the vertex $v$ (its cyclic degree is five and it has a pair of neighbors 
colored with the same color). Again, the existence of this coloring contradicts that $G$ is a counterexample.
\end{proof}

By Lemmas~\ref{lm-pentagon5} and~\ref{lm-pentagon6}, we have:

\begin{lemma}
\label{lm-pentagonS}
Let $G$ be a minimal counterexample and $v$ a pentagonal vertex
with neighbors $v_1$, $v_2$, $v_3$, $v_4$ and $v_5$ in $G$.
If the edge $v_iv_{i+1}$ is contained in a $4$-face,
then the degrees of $v_i$ and $v_{i+1}$ are at least seven.
\end{lemma}

\begin{figure}
\begin{center}
\epsfbox{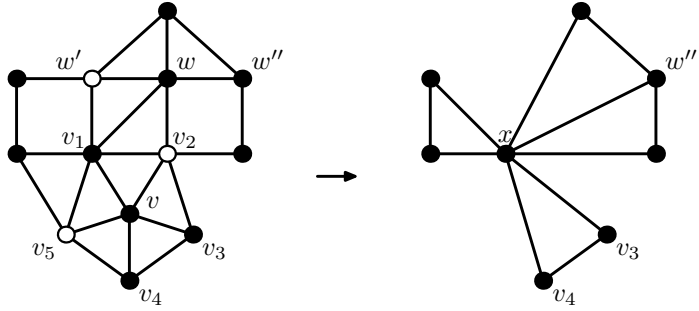}
\end{center}
\caption{The configuration described in Lemma~\ref{lm-pentagon65}. The vertex $x$ is obtained by identifying vertices drawn with empty circles.}
\label{fig-pentagon65}
\end{figure}

\begin{figure}
\begin{center}
\epsfbox{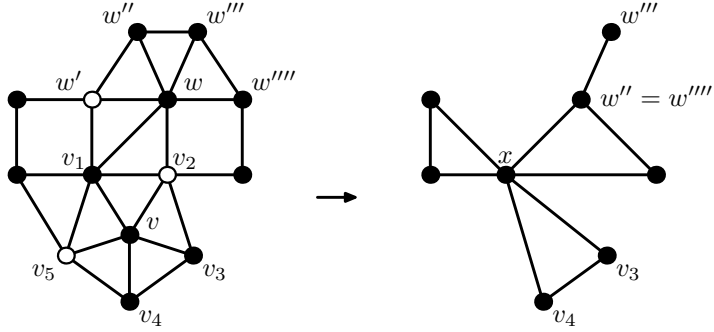}
\end{center}
\caption{The configuration described in Lemma~\ref{lm-pentagon66}. The vertex $x$ is obtained by identifying vertices drawn with empty circles.}
\label{fig-pentagon66}
\end{figure}

At the end of this section, we exclude two more complex configurations
from appearing around a pentagonal vertex in a minimal counterexample.
The configurations described in Lemmas~\ref{lm-pentagon65} and \ref{lm-pentagon66}
are depicted in Figures~\ref{fig-pentagon65} and \ref{fig-pentagon66},
respectively.

\begin{lemma}
\label{lm-pentagon65}
No minimal counterexample contains a pentagonal vertex $v$
with neighbors $v_1, \ldots, v_5$ such that for some $i\in\{1,\ldots,5\}$
\begin{enumerate}
\item the degree of $v_i$ is six,
\item the vertices $v_i$ and $v_{i+1}$ have a common neighbor $w$ of degree five,
\item the vertices $v_i$ and $w$ have a common neighbor $w'$, and
\item the edges $v_iw'$ and $v_{i+1}w$ lie in $4$-faces.
\end{enumerate}
\end{lemma}

\begin{proof}
We can assume that $i=1$. Let $w''$ be the neighbor of $w$ distinct from $v_2$
that lies on the $4$-face incident with $w$. Remove the vertices $v$, $v_1$ and
$w$ from $G$, identify the vertices $v_2$, $v_5$ and $w'$ to a new vertex $x$, and
add an edge $xw''$. Let $G'$ be the resulting graph. As any pair of identified
vertices have a common neighbor, the graph $G'\setminus\{xw''\}$ is loopless
by Lemma~\ref{lm-sep}. If the edge $xw''$ were a loop, then the vertices
$v_5$ and $w''$ would coincide in $G'$ which would yield a separating
$3$-cycle $v_1ww''=v_5$ in $G$. We conclude that $G'$ is loopless. Similarly,
all $4$-faces of $G'$ are vertex-disjoint.

By the minimality of $G$, the graph $G'$ has a cyclic $5$-coloring.
Assign vertices of $G$ the colors of their counterparts in $G'$.
The only vertices without a color are the vertices $w$, $v_1$ and $v$
which we color in this order. Let us verify that each of these vertices
is cyclically adjacent to vertices of at most four distinct colors
when we want to color it. At the beginning, the vertex $w$ has six
cyclic neighbors, out of which two have the same color ($v_2$ and $w'$) and
one is uncolored. Next, the vertex $v_1$ has cyclic degree seven but it
is adjacent to a triple of vertices with the same color and an uncolored
vertex. Finally, the cyclic degree of $v$ is five and two of its neighbors
have the same color. The constructed coloring violates our assumption that
$G$ is a counterexample.
\end{proof}

\begin{lemma}
\label{lm-pentagon66}
No minimal counterexample contains a pentagonal vertex $v$
with neighbors $v_1, \ldots, v_5$ such that for some $i\in\{1,\ldots,5\}$
\begin{enumerate}
\item the degree of $v_i$ is six,
\item the vertices $v_i$ and $v_{i+1}$ have a common neighbor $w$ of degree six,
\item the vertices $v_i$ and $w$ have a common neighbor $w'$, and
\item the edges $v_iw'$ and $v_{i+1}w$ lie in $4$-faces.
\end{enumerate}
\end{lemma}

\begin{proof}
We can assume that $i=1$. Let $w''$, $w'''$ and $w''''$ be the neighbors of $w$ 
as depicted in Figure~\ref{fig-pentagon66}.
Remove the vertices $v$, $v_1$ and
$w$ from $G$, identify the vertices $v_2$, $v_5$ and $w'$ to a new vertex $x$ and
identify the vertices $w''$ and $w''''$.
Let $G'$ be the resulting graph. As any pair of identified
vertices have a common neighbor, the graph $G'$ is loopless
by Lemma~\ref{lm-sep}. Moreover, all $4$-faces of $G'$ are vertex-disjoint.

By the minimality of $G$, the graph $G'$ has a cyclic $5$-coloring.
Now assign vertices of $G$ the colors of their counterparts in $G'$.
The only vertices without a color are the vertices $w$, $v_1$ and $v$
which we color in this order. Let us verify that each of these vertices
is cyclically adjacent to vertices of at most four distinct colors
when we want to color it. At the beginning, the vertex $w$ has seven
cyclic neighbors, out of which two pairs have the same color (the pair
$v_2$ and $w'$, and the pair $w''$ and $w''''$) and
one neighbor is uncolored. Next, the vertex $v_1$ has also cyclic degree seven
but it
is adjacent to a triple of vertices with the same color and an uncolored
vertex. Finally, the cyclic degree of $v$ is five and two of its neighbors
have the same color. Finally, the obtainec coloring contradicts that
$G$ is a counterexample.
\end{proof}

\section{Discharging rules}
\label{sect-charge}

The core of the proof is an application of the standard discharging method.
We fix a minimal counterexample and assign
each vertex and each face initial charge
as follows: each $d$-vertex receives $d-6$ units of charge and
each $d$-face receives $2d-6$ units of charge. An easy application of Euler
formula yields that the sum of initial amounts of charge is $-12$.
The amount of charge is then redistributed using the rules introduced
in this section in such a way that all vertices and faces have non-negative
amount of charge at the end. Since the redistribution preserves the total
amount of charge, this will eventually contradict the existence of
a minimal counterexample.

Let us start presenting the rules for charge redistribution.
Rules S1 and S2 guarantee that the amount of final charge of
every vertex incident with a $4$-face is zero (vertices not
incident with a $4$-face are not affected by Rules S1 and S2).

\begin{description}
\item[Rule S1] Every $5$-vertex receives $1$ unit of charge from its (unique) incident $4$-face.
\item[Rule S2] Every $d$-vertex, $d\ge 6$, sends $d-6$ units of charge to its incident $4$-face.
\end{description}

\begin{figure}
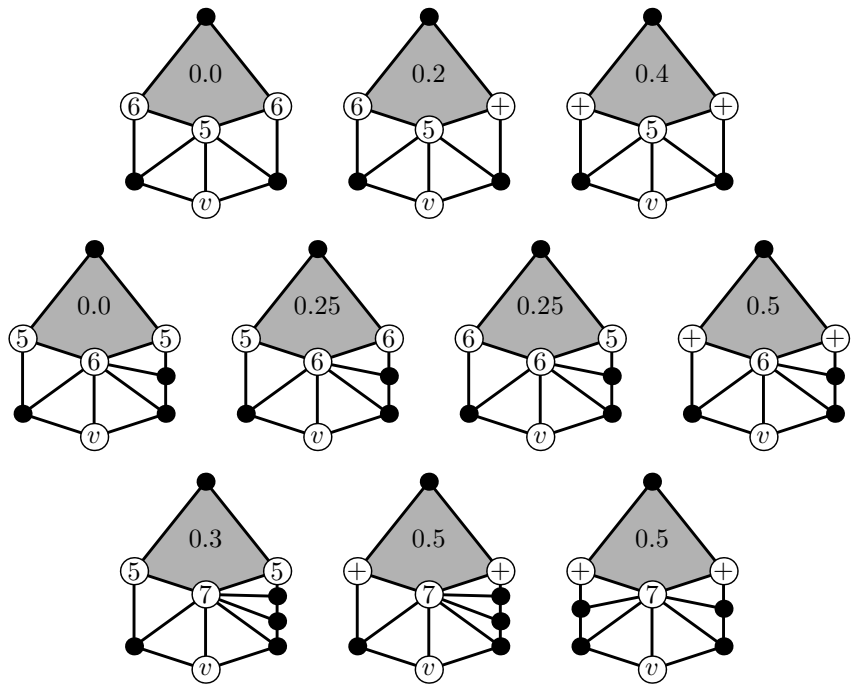

\begin{center}
\epsfbox{plan-ind.2}
\hskip 5mm
\epsfbox{plan-ind.3}
\hskip 5mm
\epsfbox{plan-ind.4}
\vskip 2mm
\epsfbox{plan-ind.5}
\hskip 5mm
\epsfbox{plan-ind.6}
\hskip 5mm
\epsfbox{plan-ind.7}
\hskip 5mm
\epsfbox{plan-ind.8}
\vskip 2mm
\epsfbox{plan-ind.9}
\hskip 5mm
\epsfbox{plan-ind.10}
\hskip 5mm
\epsfbox{plan-ind.11}
\end{center}
\caption{Illustration of Rules P5a--P7c. The numbers in circles represent
         degrees of vertices (plus signs stand for any degree not constrained
         in another part of the figure), the $4$-face $f$ sending charge is shaded and
	 the pentagonal vertex receiving charge is denoted by $v$.
	 The amount of charge sent is represented by the number
	 in the middle of the face $f$.}
\label{fig-rules}
\end{figure}

A more complex set of rules is needed to guarantee that the amount of final
charge of pentagonal vertices is non-negative. The following notation
is used in Rules P5a--P8+: $v$ is a pentagonal vertex adjacent to a vertex $w$
incident with a $4$-face $f$ distant from $v$; the neighbors of $w$
incident with $f$ are denoted $w'$ and $w''$. A vertex $w$ is understood
to be one-sided or double-sided with respect to $v$. Rules P5a--P7c are
illustrated in Figure~\ref{fig-rules}.

\begin{description}
\item[Rule PC] Every pentagonal vertex receives $1$ unit of charge
               from each close $4$-face.
\item[Rule P5a] If $w$ has degree five and
                exactly one of the vertices $w'$ and $w''$ have degree six,
                then $v$ receives $0.2$ units of charge from $f$.
\item[Rule P5b] If $w$ has degree five and both $w'$ and $w''$ have
                degree at least seven,
                then $v$ receives $0.4$ units of charge from $f$.
\item[Rule P6a] If $w$ has degree six,
                exactly one of the vertices $w'$ and $w''$ have degree five and
		the other has degree six,
                then $v$ receives $0.25$ units of charge from $f$.
\item[Rule P6b] If $w$ has degree six and
                the sum of the degrees of $w'$ and $w''$ is at least twelve,
                then $v$ receives $0.5$ units of charge from $f$.
\item[Rule P7a]	If $w$ is a one-sided vertex of degree seven and
                both $w'$ and $w''$ have degree five,
		then $v$ receives $0.3$ units of charge from $f$.
\item[Rule P7b]	If $w$ is a one-sided vertex of degree seven and
                at most one of the vertices $w'$ and $w''$ has degree five,
		then $v$ receives $0.5$ units of charge from $f$.
\item[Rule P7c] If $w$ is a double-sided vertex of degree seven,
		then $v$ receives $0.5$ units of charge from $f$.
\item[Rule P8+] If the degree of $w$ is eight or more,
               then $v$ receives $0.5$ units of charge from $f$.
\end{description}

The amount of final charge of faces and vertices after redistributing
charge based on the above rules is analyzed in the next two sections.

\section{Final charge of faces}
\label{sect-faces}

In this section, we analyze the final amount of charge of faces
in a minimal counterexample. Since $3$-faces do not receive or send out
any charge, it is enough to analyze the final charge of $4$-faces.
We break down the analysis into four lemmas that cover all possible
cases how a $4$-face can look like (up to symmetry). We start with $4$-faces
incident with two vertices of degree five.

\begin{lemma}
\label{lm-f55}
Let $f=v_1v_2v_3v_4$ be a $4$-face of a minimal counterexample.
If the degrees of $v_1$ and $v_3$ are five,
then the final amount of charge of $f$ is non-negative.
\end{lemma}

\begin{proof}
By Lemma~\ref{lm-square}, the degree of $v_2$ and $v_4$ is at least six, and
by Lemma~\ref{lm-pentagonS}, no pentagonal vertex is close to $f$.
Let $k$ be the number of vertices of degree seven or more incident with $f$.
By Rules P5a or P5b, the face $f$ sends pentagonal vertices adjacent
to $v_1$ or $v_3$ at most $2\times k\times 0.2=0.4k$ units of charge.
Let $d_i$ be the degree of a vertex $v_i$, $i=2,4$.
If $d_i=6$ for $i=2,4$, then $f$ sends out no charge to pentagonal vertices
adjacent to $v_i$. If $d_i=7$ for $i=2,4$,
then the face $f$ sends
either $0.3$ units of charge to at most two pentagonal vertices
adjacent to $v_i$ by Rule P7a or $0.5$ units of charge to
a single vertex by Rule P7c; this follows from the fact
no two adjacent neighbors of a vertex $v_i$ can be both pentagonal and
the common neighbors of $v_i$ and $v_1$ or $v_3$ are not pentagonal
by Lemma~\ref{lm-pentagon5}. These two facts also imply for $d_i>7$ that
$v_i$ sends to each of at most $(d_i-3)/2$ pentagonal vertices
adjacent to $v_i$ $0.5$ units of charge by Rule P8+.

Let us summarize. After Rules S1 and S2 apply, the amount of charge of $f$
is equal to $d_2+d_4-12$. We next distinguish several cases based on $d_2$ and
$d_4$:
\begin{itemize}
\item If $d_2=6$ and $d_4=6$, no further charge is sent
      out and the final charge of $f$ is zero.
\item If $d_2=6$ and $d_4=7$ (or vice versa),
      $f$ sends out at most $0.4$ units of charge to pentagonal vertices
      adjacent to $v_1$ or $v_3$ and at most $0.6$ units of charge to such
      vertices adjacent to $v_4$. Hence, its final charge is again non-negative.
\item If $d_2=6$ and $d_4>7$ (or vice versa),
      $f$ sends out at most $0.4$ units of charge to pentagonal vertices
      adjacent to $v_1$ or $v_3$ and at most $(d_4-3)/4$ units of charge
      to such vertices adjacent to $v_4$.
      Hence, its final charge is again non-negative.
\item If $d_2=7$ and $d_4=7$,
      $f$ sends out at most $0.8$ units of charge to pentagonal vertices
      adjacent to $v_1$ or $v_3$, at most $0.6$ units of charge to pentagonal
      vertices adjacent to $v_2$ and at most $0.6$ units of charge
      to pentagonal vertices adjacent to $v_4$.
      Its final charge is again non-negative.
\item If $d_2=7$ and $d_4>7$ (or vice versa),
      $f$ sends out at most $0.8$ units of charge to pentagonal vertices
      adjacent to $v_1$ or $v_3$, at most $0.6$ units of charge to such
      vertices adjacent to $v_2$ and at most $(d_4-3)/4$ units of charge
      to pentagonal vertices adjacent to $v_4$.
      Hence, its final charge is again non-negative.
\item If $d_2>7$ and $d_4>7$,
      the face $f$ sends out at most $0.8$ units of charge to pentagonal vertices
      adjacent to $v_1$ or $v_3$, and at most $(d_2+d_4-6)/4$ units of charge
      to such vertices adjacent to $v_2$ or $v_4$.
      Hence, its final charge is again non-negative.
\end{itemize}
\end{proof}

Next, we analyze $4$-faces incident with vertices of degree seven or more only.
Note that the bound on the number of pentagonal neighbors of vertices of a $4$-face
is also used in Lemmas~\ref{lm-f5}--\ref{lm-f6} without giving so much details
on its derivation as in the proof of Lemma~\ref{lm-f7}.

\begin{lemma}
\label{lm-f7}
Let $f=v_1v_2v_3v_4$ be a $4$-face of a minimal counterexample.
If the degrees of $v_1$, $v_2$, $v_3$ and $v_4$ are at least seven,
then the final amount of charge of $f$ is non-negative.
\end{lemma}

\begin{proof}
Let $D$ be the sum of the degrees of the vertices $v_1$, $v_2$, $v_3$ and
$v_4$. After Rule S2 applies to each of these four vertices, the face $f$
has charge $D-22$. Rules PC, P7a, P7b, P7c and P8+ apply at most $(D-12)/2$
vertices. The vertices $v_1$, $v_2$, $v_3$ and $v_4$ have $D-8$ neighbors
not incident with the face $f$ counting the common neighbors of them twice.
Hence, if the common neighbors of $v_i$ and $v_{i+1}$ are counted once,
there are at most $D-12$ neighbors not incident with $f$ and since no two
adjacent vertices can be both pentagonal, the number of pentagonal neighbors
is at most $(D-12)/2$.

Rule PC can apply at most $4$ times since
a single $4$-face can be close to at most $4$ pentagonal vertices.
Since $f$ can send out at most $0.5$ units of charge by Rules P7a, P7b, P7c and
P8+, and it can send out at most $1$ unit of charge by Rule PC,
the $4$-face $f$ sends out at most the following amount of charge:
$$\frac{D-12}{2}\times 0.5+4\times 0.50=\frac{D}{4}-1\;\mbox{.}$$
By the assumptions of the lemma, the degree of each vertex $v_i$
is at least $7$ and thus $D\ge 28$. Since $D/4-1\le D-22$ for $D\ge 28$,
the final amount of charge of $f$ is non-negative.
\end{proof}

We next analyze $4$-faces incident with a single vertex of degree five.

\begin{lemma}
\label{lm-f5}
Let $f=v_1v_2v_3v_4$ be a $4$-face of a minimal counterexample.
If the degree of $v_1$ is five and the degree of $v_3$ is at least six,
then the final amount of charge of $f$ is non-negative.
\end{lemma}

\begin{proof}
If all vertices $v_2$, $v_3$ and $v_4$ have degree six,
then $f$ can send out $0.25$ units of charge by Rule P6a
to pentagonal neighbors of $v_2$ and $v_4$ (note that each of these
two vertices has at most one such pentagonal neighbor) and
$0.5$ units of charge by Rule P6b to a pentagonal neighbor of $v_3$.
Observe that no pentagonal vertex is close to $f$
by Lemma~\ref{lm-pentagonS}.
Altogether, $f$ receives no charge and sends out at most $2$ units of
charge (one unit by Rule S1 to $v_1$). Consequently, its final charge is
non-negative.

If two of the vertices $v_2$, $v_3$ and $v_4$ have degree six and
one has degree $d\ge 7$, then $f$ can send out at most $0.2$ units
of charge to a pentagonal neighbor of $v_1$, at most $0.5$ units
charge to a pentagonal neighbor of each vertex of degree six,
at most $0.5$ to at most $(d-3)/2$ pentagonal neighbors of the vertex
of degree $d$ and $1$ unit of charge to $v_1$. Altogether, it sends
out at most $(d-3)/4+2.2=d/4+1.45$ units of charge. Since the initial
charge of $f$ amounts to $2$ units and $f$ receives $d-6$ units by Rule S2,
its final charge is non-negative if $d\ge 8$ (observe that $d/4+1.45\le d-4$
for $d\ge 8$).
If $d=7$ and the vertex of degree $d$ is $v_2$, then $f$ can send
$1$ unit of charge to $v_1$ by Rule S1, $0.2$ units of charge to
a pentagonal neighbor of $v_1$ by Rule P5b, $0.5$ units of charge to
each of at most two pentagonal neighbors of $v_2$ by Rule P7b or P7c,
$0.5$ units of charge to a pentagonal neighbor of $v_3$ by Rule P6b and
$0.25$ units of charge to a pentagonal neighbor of $v_4$ by Rule P6a.
In total, $f$ sends out at most $2.95$ units of charge. The case that
the vertex of degree $d=7$ is $v_4$ is symmetric to this one.
Finally, if the vertex of degree $d=7$ is $v_3$, then $f$ can send
$1$ unit of charge to $v_1$ by Rule S1 and $0.5$ units of charge
to at most four pentagonal neighbors of $v_2$, $v_3$ and $v_4$.
The face $f$ sends no charge to a pentagonal neighbor of $v_1$
since neither Rule P5a nor P5b can apply. Again, the final charge
of $f$ is non-negative.
%If $d=7$, then the face $f$ sends at most $0.6$ units
%of charge to pentagonal neighbors of the vertex of degree $d=7$ (either
%twice $0.3$ units by Rule P7a or $0.5$ units by Rule P7b or P7c).
%Hence, the charge sent out by $f$ is at most $0.6+2.2=2.8$
%while the initial charge of $f$ equals to $2$ units and $f$ gets
%$1$ unit of charge from the vertex of degree $d=7$ by Rule S2.

We now assume that only one of the vertices $v_2$, $v_3$ and $v_4$
have degree six and the remaining two vertices have degrees $d$ and $d'$,
$d\ge 7$ and $d'\ge 7$. The face $f$ sends out $1$ unit of charge
to $v_1$ by Rule S1, at most $0.40$ units
of charge to a pentagonal neighbor of $v_1$,
at most $0.50$ units of charge to a pentagonal neighbor of
the vertex of degree six, and at most $0.50$ units of charge
to each of at most $(d+d'-6)/2$ pentagonal neighbors of vertices
of degree $d$ and $d'$ unless Rule PC applies. Rule PC can apply
at most once by Lemma~\ref{lm-pentagonS}. Since the initial amount
charge of $f$ is $2$, $f$ receives $d+d'-12$ units by Rule S2 and
sends out at most $1+(d+d'-6)/4+0.90+0.50=(d+d')/4+0.90$ units of charge and
at most $(d+d')/4+0.50$ if Rule PC does not apply,
the final charge of $f$ is non-negative (note that $(d+d')/4+0.90\le d+d'-10$
for $d+d'\ge 15$) unless $d=d'=7$ and Rule PC also applies.
If $d=d'=7$ and Rule PC applies, Lemma~\ref{lm-pentagon6} implies that
the vertices of degree seven are adjacent. By symmetry, $v_1$ has degree five,
$v_2$ has degree six and $v_3$ and $v_4$ have degree seven.
Hence, $f$ can send $1$ unit of charge to $v_1$ by Rule S1,
$0.2$ units of charge to a pentagonal neighbor of $v_1$ by Rule P5a,
at most $0.5$ units of charge to each of at most three pentagonal
neighbors of $v_2$, $v_3$ and $v_4$ that are not close and
$1$ unit of charge to the close pentagonal neighbor by Rule PC.
We conclude that $f$ sends out at most $1+0.2+3\cdot 0.5+1=3.7$ units
of charge while it receives $2$ units of charge by Rule S2 in addition
to $2$ units of its initial charge, i.e., its final charge is non-negative.

It remains to consider the case when all the vertices $v_2$, $v_3$ and $v_4$
have degree at least seven. Let $d_i$ be the degree of the vertex $v_i$,
$i=2,3,4$. There are at most $(d_2+d_3+d_4-9)/2$ pentagonal neighbors
of the vertices $v_2$, $v_3$ and $v_4$ and Rule PC can apply at most twice.
In addition, the face $f$ can send out $0.4$ units of charge to a pentagonal
neighbor of a vertex $v_1$ and $1$ unit of charge to $v_1$ by Rule S1.
Altogether, the amount of charge sent out by $f$ is at most:
$$1.4+\frac{d_2+d_3+d_4-9}{2}\times 0.5+2\times 0.5=\frac{d_2+d_3+d_4}{4}+0.15\;\mbox{.}$$
The initial amount of charge of $f$ is $2$ units and $f$
receives $d_2+d_3+d_4-18$ units of charge by Rule S2 from the vertices
$v_2$, $v_3$ and $v_4$. Hence, if $d_2+d_3+d_4\ge 22$, then the final
charge of the face $f$ is clearly non-negative.

If $d_2+d_3+d_4=21$, then all the degrees $d_2$, $d_3$ and $d_4$ must
be equal to $7$. If the vertices $v_2$, $v_3$ and $v_4$ have six pentagonal
neighbors, then none of them is close to $f$. Hence, Rule PC never applies.
We conclude that $f$ sends out at most the following amount of charge:
$$1.4+6\times 0.5=4.4\;\mbox{.}$$
On the other hand, if there are at most five pentagonal neighbors of $v_2$,
$v_3$ and $v_4$, Rule PC can apply (at most twice). Hence, the charge sent
out by $f$ is at most:
$$1.4+5\times 0.5+2\times 0.5=4.9\;\mbox{.}$$
Since the initial amount of charge of $f$ is $2$ units and $f$ receives
$3$ units of charge from the vertices $v_2$, $v_3$ and $v_4$, its final
charge is non-negative.
\end{proof}

Finally, we analyze $4$-faces incident with vertices of degree six but no vertices of degree five.

\begin{lemma}
\label{lm-f6}
Let $f=v_1v_2v_3v_4$ be a $4$-face of a minimal counterexample.
If the degree of $v_1$ is six and the degrees of $v_2$, $v_3$ and $v_4$
are at least six,
then the final amount of charge of $f$ is non-negative.
\end{lemma}

\begin{proof}
Let $D$ be the sum of the degrees of the vertices $v_1$, $v_2$, $v_3$ and
$v_4$. After Rule S2 applies to each of these four vertices, the face $f$
has charge $D-22$.
We now distinguish several cases based on which vertices $v_i$,
$i=1,2,3,4$, have degree six:
\begin{itemize}
\item If all vertices $v_i$ have degree six, then there is no pentagonal
      vertex close to $f$ by Lemma~\ref{lm-pentagon6}.
      Hence, each $v_i$ is adjacent to at most
      one pentagonal vertex and $f$ sends $0.5$ units of charge
      by Rule P6b at most four times. This implies that the final
      amount of charge of $f$ is non-negative.
\item If three vertices $v_i$ have degree six, then there is again 
      no pentagonal vertex close to $f$
      by Lemma~\ref{lm-pentagon6}. Let $d$ be the degree of
      the vertex with degree seven or more. Such vertex is adjacent
      to at most $(d-3)/2$ pentagonal vertices and each other vertex
      to at most one pentagonal vertex. Hence, $f$ sends out at most
      $(d-3)/4+3/2=d/4+3/4$ units of charge. Since its charge after applying
      Rule S2 was $D-22=d-4$ and $d\ge 7$, its final amount of charge
      is non-negative.
\item It two vertices $v_i$ have degree six, then there is at most one
      pentagonal vertex close to $f$. The charge is sent by $f$
      to at most $(D-12)/2$ pentagonal vertices and at most once
      by Rule PC.
      Hence, the total amount of charge sent out is at most
      $$\frac{D-12}{2}\times 0.5+0.5=\frac{D}{4}-2.5\;\mbox{.}$$
      Since $D\ge 26$ and the charge of $f$ after applying Rule S2
      is at $D-22$, the final amount of charge of $f$ is non-negative.
\item If $v_1$ is the only vertex $v_i$ with degree six, the charge is sent by $f$
      to at most $(D-12)/2$ pentagonal vertices and at most twice by Rule PC.
      Hence, the total amount of charge sent out is at most
      $$\frac{D-12}{2}\times 0.5+2\times 0.5=\frac{D}{4}-2\;\mbox{.}$$
      Since $D\ge 27$ and the charge of $f$ after applying Rule S2
      is at $D-22$, the final amount of charge of $f$ is non-negative.
\end{itemize}
\end{proof}

\section{Final charge of vertices}
\label{sect-vertices}

A minimal counterexample has no vertices of degree four or less
by Lemma~\ref{lm-mindeg}.
The amount of final charge of vertices that are not pentagonal
is non-negative: vertices incident with a $4$-face have zero final charge 
since only Rule S1 or S2 can apply to them and
other non-pentagonal vertices keep their original (non-negative) charge
since none of the rules applies to them (note that every vertex of degree
five is either pentagonal or incident with a $4$-face by Lemma~\ref{lm-pentagonal}).

Hence, we can focus on the amount of final charge of pentagonal
vertices. Pentagonal vertices that are not solitary receive $1$ unit
of charge from a close $4$-face by Rule PC and thus their final charge
is non-negative. We now analyze the amount of charge
of solitary pentagonal vertices and start with those adjacent
to a vertex of degree five.

\begin{lemma}
\label{lm-v5}
Every solitary pentagonal vertex $v$ adjacent to a vertex of degree five
has non-negative final charge.
\end{lemma}

\begin{figure}
\begin{center}
\epsfbox{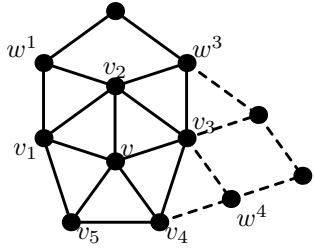}
\end{center}
\caption{Notation used in the proof of Lemma~\ref{lm-v5}.}
\label{fig-v5}
\end{figure}

\begin{proof}
Let $v_1,\ldots,v_5$ be the neighbors of $v$ and $f_i$ the $4$-face
containing the vertex $v_i$, $i=1,\ldots,5$.
By symmetry, we can assume that the degree of $v_2$ is five.
Since no two $4$-faces share a vertex, $v$ has a double-sided neighbor $v_k$.
Note that $k\not=2$ and the $4$-face $f_k$ sends
$0.5$ units of charge to $v$ (either by Rule P7c or Rule P8+).

Let $w^1$ be the common neighbor of $v_1$ and $v_2$ and $w^3$ the common
neighbor of $v_2$ and $v_3$ (see Figure~\ref{fig-v5}). Since the degree of $v_2$ is five,
the degrees of $w^1$ and $w^3$ are at least six by Lemma~\ref{lm-square}.
If the degree of $w^1$ is six, then the degree of $v_1$
is at least seven by Lemma~\ref{lm-pentagon66} and
the $4$-face $f_1$ sends $v$ at least $0.3$ units of charge.
Similarly, if the degree of $w^3$ is six, then the $4$-face
$f_3$ sends $v$ at least $0.3$ units of charge. On the other hand,
if the degree of at least one of the vertices $w^1$ and $w^3$ is bigger
than six, then $v$ receives at least $0.2$ units of charge from the $4$-face
$f_2$, and if the degrees of both $w^1$ and $w^3$ are bigger
than six, then $v$ receives at least $0.4$ units of charge from $f_2$.

We conclude that if $k\not\in\{1,3\}$, then $v$ receives $0.5$ units of charge
from $f_k$ and at least $0.4$ units of charge
from the faces $f_1$, $f_2$ and $f_3$. In particular,
the final charge of $v$ is non-negative unless $v$ receives
exactly $0.4$ units of charge from the faces $f_1$, $f_2$ and $f_3$ altogether.
In such case, $v$ receives $0.4$ units of charge from $f_2$, which implies
that the degrees of $w^1$ and $w^3$ are more than six, and no charge is sent from $f_1$ or $f_3$,
which implies that the degrees of $v_1$ and $v_3$ are six and the degrees
of their neighbors on $f_1$ and $f_3$ are five.
Let us analyze this case in more detail. By symmetry, we can assume that
$k=5$. Let $w^4$ be the common neighbor of $v_3$ and $v_4$. Since $f_3$ sends no charge,
the degree of $w^4$ is five. Hence, the degree of $v_4$
is at least seven by Lemma~\ref{lm-pentagon65}. Consequently, the face $f_4$
sends $v$ at least $0.3$ units of charge. Altogether, $v$ receives
$0.4$ units of charge from $f_2$, at least $0.3$ units of charge from $f_4$ and
$0.5$ units of charge from $f_5$ and its final charge is non-negative.
We have just shown that if $k\not\in\{1,3\}$, then the final charge
of $v$ is non-negative.

In the rest, we assume that $k=1$ and one of the following
two cases applies (otherwise, the faces $f_2$ and $f_3$ send together
at least $0.5$ units of charge to $v$ and the final charge of $v$ is
non-negative). The other cases are excluded by Lemmas~\ref{lm-square} and
\ref{lm-pentagon66}.
\begin{itemize}
\item {\bf The $4$-face $f_2$ sends $v$ no charge and
      the $4$-face $f_3$ sends $0.3$ units of charge.}\\
      In this case, the degrees of both $w^1$ and $w^3$ are six and
      $v_3$ is a one-sided vertex with degree seven with both neighbors
      on $f_3$ of degree five. In particular,
      the common neighbor $w^4$ of $v_3$ and $v_4$ lies in the face $f_3$ and
      it has degree five.
      By Lemma~\ref{lm-pentagon65}, the degree of $v_4$ is at least seven and
      thus the $4$-face $f_4$ sends at least $0.3$ units
      of charge to $v$. In total, $v$ receives $0.5$ units of charge
      from $f_1$, $0.3$ units of charge from $f_3$ and
      at least $0.3$ units of charge from $f_4$. We conclude that
      the final charge of $v$ is non-negative.
\item {\bf The $4$-face $f_2$ sends $0.2$ or $0.4$ units of charge and
      the $4$-face $f_3$ sends no charge.}\\
      In this case, $v_3$ has degree six and its common neighbor $w^4$
      with the vertex $v_4$ has degree five and lies on the face $f_3$.
      Lemma~\ref{lm-pentagon65} now implies that the degree of $v_4$
      is at least seven. Hence, the face $f_4$ sends at least $0.3$ units
      of charge to $v$. Summarizing, $v$ receives $0.5$ units of charge
      from $f_1$, at least $0.2$ units of charge from $f_2$ and
      at least $0.3$ units of charge from $f_4$ which makes its final
      charge non-negative.
\item {\bf The $4$-face $f_2$ sends $0.2$ units of charge and
      the $4$-face $f_3$ sends $0.25$ units of charge.}\\
      In this case, $v_3$ has degree six and its common neighbor $w^4$
      with the vertex $v_4$ has degree five or six and lies on the face $f_3$.
      Lemmas~\ref{lm-pentagon65} and~\ref{lm-pentagon66} yield that
      the degree of $v_4$ is at least seven. This implies that
      the face $f_4$ sends at least $0.3$ units of charge to $v$.
      We conclude that $v$ receives $0.5$ units of charge
      from $f_1$, $0.2$ units of charge from $f_2$, $0.25$ units of charge
      from $f_3$ and at least $0.3$ units of charge from $f_4$, and
      the final charge of $v$ is non-negative.
\end{itemize}
\end{proof}

It remains to analyze solitary pentagonal vertices
adjacent to no vertices of degree five.

\begin{lemma}
\label{lm-v6}
Every solitary pentagonal vertex $v$ adjacent to no vertex of degree five
has non-negative final charge.
\end{lemma}

\begin{proof}
Let $v_1,\ldots,v_5$ be the neighbors of $v$ and $f_1,\ldots,f_5$
the $4$-faces incident with the neighbors of $v$ as in the proof
of Lemma~\ref{lm-v5}. If $v$ receives charge from at least four
of the faces $f_1,\ldots,f_5$, then it receives at least $1$ unit
of charge in total and its final charge is non-negative. Hence,
we can assume that $v$ does not receive charge from two of the faces,
by symmetry, from the face $f_1$ and the face $f_2$ or $f_3$.
Note that if $v$ receives no charge from the face $f_i$,
then $v_i$ has degree six and both its neighbors on $f_i$
must have degree five.

\begin{figure}
\begin{center}
\epsfbox{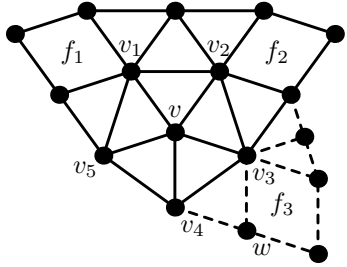}
\end{center}
\caption{Notation used in the first part of the proof of Lemma~\ref{lm-v6}.}
\label{fig-v6a}
\end{figure}

Let us first assume that
the vertex $v$ receives no charge from the faces $f_1$ and $f_2$.
The situation is depicted in Figure~\ref{fig-v6a}; note that
the vertices $v_1$ and $v_2$ cannot have a common neighbor
of degree five on a face $f_1$ or $f_2$ by Lemma~\ref{lm-pentagon65}.
Observe that there must be a double-sided vertex $v_k$, $k\in\{3,4,5\}$.
By Lemma~\ref{lm-pentagon65}, the degrees of the vertices $v_3$ and $v_5$
are at least seven. Hence, if $k=4$, $v$ receives at least $0.3$ units
of charge from the faces $f_3$ and $f_5$ and $0.5$ units of charge
from $f_4$, and its final charge is non-negative.

We now assume that $k=5$ and the face $f_3$ sends only $0.3$ units
of charge to $v$ (otherwise, $v$ receives $0.5$ units of charge
from $f_3$ and its final charge is non-negative).
Hence, $v_3$ is a one-sided vertex of degree seven and
the common neighbor $w$ of $v_3$ and $v_4$ has degree five and lies on $f_3$.
Consequently, the degree of $v_4$ is at least seven
by Lemma~\ref{lm-pentagon65}. We conclude that $v$ receives
$0.3$ units of charge from $f_3$, at least $0.3$ units of charge
from $f_4$ and $0.5$ units of charge from $f_5$. Again, the final
charge of $v$ is non-negative.

We have ruled out the case that there would be two adjacent neighbors of $v$
whose $4$-faces sent no charge to $v$. Hence, it remains to analyze the case
when the faces $f_1$ and $f_3$ send no charge to $v$.
We claim that the face $f_2$ sends $0.5$ units of charge to $v$.
This clearly holds if $v_2$ is double-sided or its degree is at least eight.
If the degree of $v_2$ is six, then $f_2$ sends $0.5$ units
of charge unless the neighbors of $v_2$ on $f_2$ have degrees five and six.
Such configurations
are excluded by Lemmas~\ref{lm-pentagon65} and~\ref{lm-pentagon66}.
Finally, if $v_2$ is one-sided and its degree is seven,
then $f_2$ sends $0.5$ units of charge to $v$
unless both the neighbors of $v_2$ on $f_2$ have degree five. One of these
neighbors is also a neighbor of $v_1$ or $v_3$ which is impossible
by Lemma~\ref{lm-pentagon66}.

We have shown that $v$ receives $0.5$ units of charge from $f_2$.
Since $v$ receives in addition at least $0.25$ units of charge from each
of the faces $f_4$ and $f_5$, its final charge is non-negative.
\end{proof}

Lemmas~\ref{lm-f55}--\ref{lm-v6} now yield Theorem~\ref{thm-main} as explained
in Section~\ref{sect-notation}.

\section{Final remarks}

If $G$ is a plane graph with faces of size three only,
then Four Color Theorem implies that $G$ is cyclically
$4$-colorable. Our theorem asserts that every plane graph
with faces of size three and four such that all faces
of size four are vertex-disjoint is cyclically $5$-colorable.
It is natural to ask whether the following might be true:

\begin{problem}
Every plane graph $G$ with maximum face size $\Delta^*$ such that
all faces of size four or more are vertex-disjoint is cyclically
$(\Delta^*+1)$-colorable.
\end{problem}

Let us remark that it is quite easy to see that such graphs $G$
are $(\Delta^*+3)$-colorable. Indeed, adding a clique to every
face of size four or more results in a graph with average
degree less than $\Delta^*+3$. After removing a vertex from $G$
that has degree less than $\Delta^*+3$ in the modified graph and
adding edges to $G$ in such a way that big faces are still
vertex-disjoint and all vertices lying on a common face in $G$
lie on a common face in the new graph, induction can be applied
to the new graph which yields the proof of the claimed bound.


\begin{thebibliography}{99}
\bibitem{bib-albertson08}
M. Albertson:
Chromatic number, independence ratio, and crossing number,
Ars Math. Contemporanea 1 (2008), 1--6.
\bibitem{bib-albertson}
M. Albertson:
Colorings and crossings,
presentation at SIAM Conference on Discrete Mathematics 2008, Burlington, VT.
\bibitem{bib-amini08+}
O. Amini, L. Esperet, J. van den Heuvel:
A unified approach to distance-two colouring of planar graphs,
manuscript.
\bibitem{bib-appel76+}
K. Appel, W. Haken:
Every planar map is four colorable,
Bull. Am. Math. Soc. 82 (1976), 449--456.
\bibitem{bib-borodin84}
O. Borodin:
Solution of Ringel's problems on vertex-face coloring of plane graphs and coloring of $1$-planar graphs,
Met. Discret. Anal. Novosibirsk 41 (1984), 12--26 (in Russian).
\bibitem{bib-borodin92}
O. Borodin:
Cyclic coloring of plane graphs,
Discrete Math. 100 (1992), 281--289.
\bibitem{bib-borodin95}
O. Borodin:
A new proof of the 6 Color Theorem,
J. Graph Theory 19 (1995), 507--521.
\bibitem{bib-borodin99+}
O. Borodin, D. P. Sanders, Y. Zhao:
On cyclic colorings and their generalizations,
Discrete Math. 203 (1999), 23--40.
\bibitem{bib-enomoto01+}
H. Enomoto, M. Hor{\v n}{\'a}k, S. Jendrol':
Cyclic chromatic number of 3-connected plane graphs,
SIAM. J. Discrete Math. 14 (2001), 121--137.
\bibitem{bib-havet07+}
F. Havet, J. van den Heuvel, C. McDiarmid, B. Reed:
List colouring squares of planar graphs, in:
Electronic Notes in Discrete Mathematics  29 (2007), 515--519
(Proceedings of EuroComb'07).
\bibitem{bib-havet08+}
F. Havet, J. van den Heuvel, C. McDiarmid, B. Reed:
List colouring squares of planar graphs,
manuscript.
\bibitem{bib-havet+}
F. Havet, D. Kr{\'a}l', J.-S. Sereni, R. {\v S}krekovski:
Facial colorings using Hall's Theorem,
submitted.
\bibitem{bib-havet++}
F. Havet, J.-S. Sereni, R. {\v S}krekovski:
$3$-facial colouring of plane graphs,
SIAM. J. Discrete Math. 22 (2008), 231--247.
\bibitem{bib-hornak99+}
M. Hor{\v n}{\'a}k, S. Jendrol':
On a conjecture by Plummer and Toft,
J. Graph Theory 30 (1999), 177--189.
\bibitem{bib-hornak00+}
M. Hor{\v n}{\'a}k, S. Jendrol':
On vertex types and cyclic colourings of 3-connected plane graphs,
Discrete Math. 212 (2000), 101--109.
\bibitem{bib-hornak+}
M. Hor{\v n}{\'a}k, J. Zl{\'a}malov{\'a}:
Another step towards proving a conjecture of Plummer and Toft,
submitted.
\bibitem{bib-kral05+}
D. Kr{\'a}l', T. Madaras, R. {\v S}krekovski:
Cyclic, diagonal and facial coloring,
European J. Combin. 26 (2005), 473--490.
\bibitem{bib-kral07+}
D. Kr{\'a}l', T. Madaras, R. {\v S}krekovski:
Cyclic, diagonal and facial coloring---a missing case,
European J. Combin. 28 (2007), 1637--1639.
\bibitem{bib-ore69+}
O. Ore, M. D. Plummer:
Cyclic coloration of plane graphs,
in: Recent progress in combinatorics (Proceedings of the Third Waterloo Conference on Combinatorics, May 1968) (W. T. Tutte, ed.), Academic Press, 1969.
\bibitem{bib-plummer87+}
M. D. Plummer, B. Toft:
Cyclic coloration of 3-polytopes,
J. Graph Theory 11 (1987), 507--515.
\bibitem{bib-ringel65}
G. Ringel:
Ein Sechsfarbenproblem auf der Kugel,
Abh. Math. Sem. Univ. Hamburg 29 (1965), 107--117 (in German).
\bibitem{bib-robertson97+}
N. Robertson, D. Sanders, D. Seymour, R. Thomas:
The four color theorem,
J. Combin. Theory Ser. B 70 (1997), 2--44.
\bibitem{bib-sanders01+}
D. P. Sanders, Y. Zhao:
A new bound on the cyclic chromatic number,
J. Combin. Theory Ser. B 83 (2001), 102--111.
\end{thebibliography}
\end{document}